\input amstex
\input amsppt.sty
\magnification=1200

\TagsOnRight
\NoRunningHeads  
\NoBlackBoxes

\define\Y{\Bbb Y}
\define\R{\Bbb R}
\define\Z{\Bbb Z}
\define\C{\Bbb C}

\define\al{\alpha}
\define\be{\beta}
\define\la{\lambda}

\def\Mz{M_{z,z'}}
\define\Mnz{M^{(n)}_{z,z'}}
\define\Mzxi{M_{z,z',\xi}}
\define\Pzxi{\Cal P_{z,z',\xi}}
\define\pitxi{\pi_{t,\xi}}
\define\Conf{\operatorname{Conf}}

\define\w{\tfrac{\xi}{\xi-1}}

\define\Ga{\Gamma}

\define\wt{\widetilde}

\define\sgn{\operatorname{sgn}}

\topmatter
\title
Z--Measures on partitions, Robinson--Schensted--Knuth correspondence, and $\beta=2$ random matrix ensembles
\endtitle
\date Preliminary version. May 18, 1999\enddate
\author Alexei Borodin and Grigori Olshanski
\endauthor
\abstract We suggest an hierarchy of all the results known so far about the connection of the asymptotics of combinatorial or representation theoretic problems with ``$\beta=2$ ensembles'' arising in the random matrix theory. We show that all such results are, essentially, degenerations of one general situation arising from so--called generalized regular representations of the infinite symmetric group.
\endabstract
\date Preliminary version \enddate
\endtopmatter

\document
\head Introduction \endhead
 
In last few years there appeared a number of papers indicating a strong connection of certain asymptotic problems of enumerative combinatorics and representation theory of symmetric groups with the random matrix theory, see \cite{BDJ1}, \cite{BDJ2}, \cite{BR1}, \cite{BR2},  \cite{B}, \cite{BO1}, \cite{BO2}, \cite{BOO}, \cite{P.I--V}, \cite{J1}, \cite{J2}, \cite{O}, \cite{TW3}, \cite{TW4}, and the list is definitely not complete. Such connection was also anticipated in earlier works \cite{Re}, \cite{K1}, \cite{K2}.

In this paper we suggest an hierarchy of all the results known so far about the connection of the asymptotics of combinatorial or representation theoretic problems with so-called ``$\beta=2$ ensembles'' arising in the random matrix theory\footnote{These ensembles are characterized by the property that their correlation functions have determinantal form with a scalar kernel, see below.}. We show that all such results are, essentially, degenerations of one general situation arising from so--called generalized regular representations of the infinite symmetric group, see \cite{KOV} and \S3 below.

 It is worth noting that though the hierarchy 
provides a clear understanding why this or that problem should have this or that asymptotics, the technical aspects of the proof are sometimes nontrivial and should not be underestimated.

Many claims cited below were recently proved by Kurt Johansson, we would like to thank him for keeping us informed about his work.

\head \S1. Z--Measures \endhead

For $n=1,2,\dots$, let $\Y_n$ denote the set of partitions of $n$,
which will be identified with Young diagrams with $n$ boxes. We agree
that $\Y_0$ consists of a single element --- the zero partition or
the empty diagram $\varnothing$. 

Given $\la\in\Y_n$, we write $|\la|=n$ and denote by $d=d(\la)$ the
number of diagonal boxes in $\la$. We shall use the Frobenius
notation \cite{Ma} 
$$
\la=(p_1,\dots,p_d\,|\,q_1,\dots,q_d).
$$
Here $p_i=\la_i-i$ is the number of boxes in the $i$th row of $\la$ on the
right of the $i$th diagonal box; likewise, $q_i=\la'_i-i$ is the number of 
boxes in the $i$th column of $\la$ below the $i$th diagonal box ($\la'$
stands for the transposed diagram). 

Note that 
$$
p_1>\dots>p_d\ge0, \qquad
q_1>\dots>q_d\ge0, \qquad
\sum_{i=1}^d(p_i+q_i+1)=|\la|.
$$
The numbers $p_i$, $q_i$ are called the {\it Frobenius coordinates\/}
of the diagram $\la$. 

Let $b=(i,j)$ be a box of $\la$; here $i,j$ are the row number and the column number of $b$. Recall the definition of the {\it content} and the {\it hook length} of $b$:
$$
c(b)=j-i,\quad h(b)=(\la_i-j)+(\la'_j-i)+1.
$$

We will consider two complex parameters $z,z'$ such that
the numbers $(z)_k(z')_k$ and $(-z)_k(-z')_k$ are real and strictly positive for any $k=1,2,\dots$. Here and below 
$$
(a)_k=a(a+1)\dots(a+k-1), \qquad (a)_0=1, 
$$ 
denotes the Pochhammer symbol. 

The above assumption on $z,z'$ means that one of the
following two conditions holds:

$\bullet$ either $z'=\bar z$ and $z\in\C\setminus\Z$

$\bullet$ or $z,z'\in\R$ and there exists $m\in\Z$ such that
$m<z,z'<m+1$  

We set 
$$
t=zz'
$$
and note that $t>0$.

For a Young diagram $\la$ let $\dim\la$ denote the number of the
standard Young tableaux of shape $\la$. Equivalently, $\dim\la$ is
the dimension of the irreducible representation (of the symmetric
group of degree $|\la|$) corresponding to $\la$,  see
\cite{Ma}. The well--known {\it hook formula} for $\dim\la$
has the form, see, e.g., \cite{Ma},
$$
\dim\la=\frac{|\la|!}{\prod_{b\in\la}h(b)}.
$$

We introduce a function on the Young diagrams depending on the
parameters $z,z'$:
$$
\Mz(\la)=\frac{\prod_{b\in\la}(c(b)+z)(c(b)+z')}{(t)_{|\la|}}\cdot\frac{\dim^2\la}{|\la|!}=\frac{|\la|!}{(t)_{|\la|}}\cdot\prod_{b\in\la}\frac{(c(b)+z)(c(b)+z')}{h^2(b)}.
\tag 1.1
$$

We agree that $\Mz(\varnothing)=1$. Thanks to our assumption on the
parameters, $\Mz(\la)>0$ for all $\la$. 

\proclaim{Proposition 1.1} For any $n$,
$$
\sum_{\la\in\Y_n}\Mz(\la)=1,
$$
so that the restriction of $\Mz$ to $\Y_n$ is a probability
distribution on $\Y_n$.
\endproclaim

We shall denote this distribution by $\Mnz$ and call it the {\it $n$th level z--measure}.

Proposition 1.1 is an easy corollary of Proposition 3.1 below.

Let $\Y=\Y_0\sqcup\Y_1\sqcup\dots$ denote the set of all Young
diagrams. Consider the negative binomial distribution on the
nonnegative integers, which depends on $t$ and the additional parameter
$\xi$, $0<\xi<1$:
$$
\pitxi(n)=(1-\xi)^t\,\frac{(t)_n}{n!}\,\xi^n\,,
\qquad n=0,1,\dots\,.
$$
For $\la\in\Y$ we set
$$
\Mzxi(\la)=\Mz(\la)\,\pitxi(|\la|).
$$
By the construction, $\Mzxi(\cdot)$ is a probability distribution on
$\Y$, which can be viewed as a mixture of the finite distributions
$\Mnz$. {}From the formulas for $\Mz$ and $\pitxi$ we get an explicit
expression for $\Mzxi$:
$$
\gathered
\Mzxi(\la)=(1-\xi)^t\,
{\xi\,}^{|\la|}\prod_{b\in\la}\frac{(c(b)+z)(c(b)+z')}{h^2(b)}
=(1-\xi)^t\,
{\xi\,}^{\sum\limits_{i=1}^d(p_i+q_i+1)}\,
t^d\,\\
\times\prod_{i=1}^d\frac{(z+1)_{p_i}(z'+1)_{p_i}(-z+1)_{q_i}(-z'+1)_{q_i}}
{p_i!p_i!q_i!q_i!}\,
{\det}^2\left[\frac1{p_i+q_j+1}\right]\,.
\endgathered
$$

We shall call $\Mzxi$ the {\it mixed z--measure}.
Following a certain analogy with models of statistical physics (cf.
\cite{V}) one may call $(\Y,\Mzxi)$ the {\it grand canonical
ensemble.}

Let $\Z'$ denote the set of half--integers,
$$
\Z'=\Z+\tfrac12=\{\dots,-\tfrac32,-\tfrac12,\tfrac12,\tfrac32,\dots\},
$$
and let $\Z'_+$ and $\Z'_-$ be the subsets of positive and negative
half--integers, respectively.
It will be sometimes convenient to identify both
$\Z'_+$ and $\Z'_-$ with $\Z_+=\{0,1,2,\dots\}$ by making use of the
correspondence $\pm(k+\tfrac12)\leftrightarrow k$, where $k\in\Z_+$.

Denote by $\Conf(\Z')$ the space of all finite subsets of $\Z'$ which
will be called {\it configurations}. We define an embedding
$\la\mapsto X$ of the set $\Y$ of Young 
diagrams into the set $\Conf(\Z')$ of configurations in $\Z'$ as follows: 
$$
\la=(p_1,\dots,p_d\,|\,q_1,\dots,q_d)\,\mapsto\,
X=\{p_1+\tfrac12,\dots,p_d+\tfrac12,-q_1-\tfrac12,\dots,-q_d-\tfrac12\}.
\tag1.2
$$
Under the identification $\Z'\simeq \Z_+\sqcup \Z_+$, the map
$\la\mapsto X$ is simply associating to $\la$ the collection of its
Frobenius coordinates. The image of the map consists exactly of the
configurations $X$ with the property $|X\cap\Z'_+ |=|X\cap \Z'_-|$. We call such configurations {\it balanced}.

 Under the embedding $\la\mapsto X$ the probability measure $\Mzxi$ on
$\Y$ turns into a probability
measure on the balanced configurations in $\Z'$. Following the conventional 
terminology, see \cite{DVJ}, we get a point process on $\Z'$; let us
denote it as $\Pzxi$.

The {\it $n$th correlation function} $\rho_n^{(z,z'\xi)}(x_1,\dots,x_n)$ of $\Pzxi$ is the probability that the random point configuration contains the points $\{x_1,\dots,x_n\}$. 

In \cite{BO2} we have computed all the correlation functions of $\Pzxi$. To state the result we need some notation. 

Consider the following functions in $u$ depending on $z$,
$z'$, $\xi$ as parameters, cf. \cite{BO2}, 
$$
\psi_\pm(u)=t^{1/2}\,\xi^{u+1/2}\,(1-\xi)^{\pm(z+z')}\,
\frac{\Ga(u+1\pm z)\Ga(u+1\pm z')}
{\Ga(1\pm z)\Ga(1\pm z')\Ga(u+1)\Ga(u+1)};   
$$
$$
P_\pm(u)=(\psi_\pm(u))^{1/2}\, F(\mp z,\mp z';u+1;\w), 
$$
$$
Q_\pm(u)=\frac{t^{1/2}\xi^{1/2}\,(\psi_\pm(u))^{1/2}}{1-\xi}\,
\frac{F(1\mp z,1\mp z';u+2;\w)}{u+1}.
$$

Here $F(a,b;c;w)$ is the Gauss hypergeometric function.

\proclaim{Theorem 1.2 (\cite{BO2})} The correlation functions of $\Pzxi$ have the form
$$
\rho_n^{(z,z'\xi)}(x_1,\dots,x_n)=\det[K(x_i,x_j)]_{i,j=1}^n,\qquad x_1,\dots,x_n\in\Z',
$$ 
where 
$$
K(x,y)=\frac{F_1(x)G_1(y)+F_2(x)G_2(y)}{x-y},
\tag 1.3
$$
with
$$
\gathered 
F_1(x)=\cases-{Q_+(x-\frac 12)} ,&x>0\\  
{P_-(-x-\frac 12)},&x<0\endcases;\qquad
F_2(x)=\cases {P_+(x-\frac 12)} ,&x>0\\ 
{Q_-(-x-\frac 12)} ,&x<0\endcases;
\\
G_1(x)=\cases {P_+(x-\frac 12)},&x>0\\  
-{Q_-(-x-\frac 12)},&x<0\endcases;\qquad
G_2(x)=\cases {Q_+(x-\frac 12)} ,&x>0\\  
{P_-(-x-\frac 12)},&x<0\endcases.
\endgathered
\tag 1.4
$$
\endproclaim

We call $K(x,y)$ the {\it hypergeometric kernel}.

\example{Remark 1.3}

1. The hypergeometric kernel has no singularity on the diagonal: the numerator of (1.3) vanishes if $x=y$.

2. The hypergeometric kernel satisfies the relation
$$
K(x,y)=\sgn(x)\sgn(y)K(y,x).
\tag 1.5
$$
This shows that the kernel is Hermitian with respect to the indefinite inner product in $\ell^2(\Z')=\ell^2(\Z'_+)\oplus \ell^2(\Z'_-)$ given by the operator $\operatorname{id}\oplus\operatorname{(--id)}$.

3. The restriction of the hypergeometric kernel to $\Z'_+$ has the form
$$
\frac{P_+(x-\frac 12)Q_+(y-\frac 12)-P_+(y-\frac 12)Q_+(x-\frac 12)}{x-y}.
$$
Note that this kernel is symmetric. We will call it the {\it  positive part} of the hypergeometric kernel.

4. Kernels with the symmetry (1.5) appeared before in works of mathematical physicists on solvable models of systems with positive and negative charged particles, see \cite{AF}, \cite{CJ1}, \cite{CJ2}, \cite{G}, \cite{F2}--\cite{F4} and references therein. The mixed z--measure can also be interpreted as a model for positive and negative particles on $\Z'$: positive particles may occupy locations in $\Z'_+$, negative --- in $\Z'_-$. The square of the Cauchy determinant
$$
{\det}^2\left[\frac1{p_i+q_j+1}\right]=\frac{\prod_{i<j}[(p_i-p_j)(q_i-q_j)]^2}{\prod_{i,j}(p_i+q_j+1)^2}
$$ 
in the formula for $\Mzxi$ above encodes the logarithmic interaction of the charged particles. 
\endexample

\head \S2. Three versions of Robinson--Schensted--Knuth correspondence \endhead

A description of the RSK--algorithm can be found in \cite{Fu}, \cite{Sa}.  

We start with the ``widest'' version of the RSK--correspondence due to Knuth \cite{Kn}.

Denote by $B_{k,l}^n$ the set of bijections between two sets of size $n$, the first set consists of (possibly repeated) numbers from $1$ to $k$ and the second set consists of (possibly repeated) numbers from $1$ to $l$. Such bijections are in one--to--one correspondence with matrices of size $k\times l$ with nonnegative integral entries, total sum of  entries equal to $n$: the $(i,j)$--entry shows how many times the element $i\in\{1,\dots,k\}$ is associated with the element $j\in\{1,\dots,l\}$. Clearly, $|B_{k,l}^n|=\binom {kl+n}{n}$.
  
The RSK--algorithm establishes a bijection of $B_{k,l}^n$ and the set of ordered pairs of semi--standard Young tableaux\footnote{Recall that the term `semi--standard Young tableau' stands for a tableau whose entries are weakly increasing along the rows and strictly increasing along the columns. In a standard tableau we have strict increasing in both directions.} of the same shape with $n$ boxes, the first tableau has entries from the set $\{1,\dots,k\}$, while the second --- from the set $\{1,\dots,l\}$.

 As is well--known, the number of semi--standard Young tableaux of shape $\lambda$ with entries from $\{1,\dots,k\}$ is equal to the value of the Schur symmetric function $s_\lambda(1,1,\dots,1,0,0,\dots)$ where the number of 1's equals $k$. This value can be written in the following form, see, e.g., \cite{Ma, I.3, Ex. 4},
$$
s_\lambda(\underbrace{1,1,\dots,1}_k,0,0,\dots)=\prod_{b\in\la}\frac{c(b)+k}{h(b)}.
$$
Recall also that the number of standard Young tableaux of shape $\la$ is $\dim\la$.

Hence, if we consider the uniform probability distribution on $B_{k,l}^n$, then, with respect to its image on the set of Young diagrams with $n$ boxes, the probability of a Young diagram $\la\in \Y_n$ equals
$$
\binom {kl+n}{n}^{-1}\prod_{b\in\la}\frac{(c(b)+k)(c(b)+l)}{h^2(b)}.
$$
Comparing this with (1.1) we conclude that this distribution coincides with $\Mnz$ for $z=k$, $z'=l$.

Note that these values of $z,z'$ do not satisfy our conditions on the parameters imposed in \S1. The reason is that for such $z,z'$ the values of $\Mnz$ can be zero, for example $\Mnz(\la)=0$ for all $\la$ with length (number of nonzero parts) greater than $\min\{k,l\}$. However, all values of $\Mnz$ remain nonnegative. We consider such situation as a specific degeneration of the regular picture (when the values $\Mnz$ are strictly positive).

Two other (earlier) version of the RSK--correspondence are due to Robinson and Schensted \cite{Ro}, \cite{S}.

Denote by $B_{k,\infty}^n$ the set of words of length $n$ built from the alphabet $\{1,\dots,k\}$ (our notation will become clear soon). It is a subset of $B_{k,n}^n$ characterized by the property that the numbers in the second set are all distinct (they encode the order of letters $\{1,\dots,k\}$ in the word). It means that in the corresponding matrices of size $k\times n$ every column has exactly one nonzero element which is equal to 1.  Obviously, $|B_{k,\infty}^n|=k^n$.

In this case the RSK--algorithm establishes a bijection of $B_{k,\infty}^n$ and the set of ordered pairs of Young tableaux of the same shape with $n$ boxes; the first tableau is semi--standard and it is filled with numbers from 1 to $k$, and the second tableau is standard. This means that the probability of a Young diagram $\la\in \Y_n$ with respect to the image of the uniform distribution on $B_{k,\infty}^n$ equals
$$
k^{-n}\prod_{b\in\la}\frac{c(b)+k}{h(b)}\cdot\dim\la.
$$
It is easy to see from (1.1) that this is the limit of $\Mnz$
for $z=k$ and $z'\to\infty$.

Finally, if we forbid for both sets in the definition of $B_{k,l}^n$ to have repetitions, then we get the symmetric group $S_n$. It would be logical to denote the symmetric group by $B_{\infty,\infty}^n$, see below. In the language of matrices, it means that we consider $n\times n$ matrices with 0's and 1's such that in each row and each column there is exactly one nonzero element. Clearly, $|S_n|=n!$.

The RSK--algorithm provides a bijection of the set of permutations of $n$ symbols and the set of ordered pairs of standard Young tableaux of the same shape with $n$ boxes. 
Hence, the probability of a Young diagram $\la\in \Y_n$ with respect to the distribution coming from the uniform distribution on $S_n$ equals ${\dim^2\la}/{n!}$. This distribution on the Young diagrams is called the {\it Plancherel distribution}. The relation (1.1) easily implies that the Plancherel distribution is the limit of $\Mnz$
when $z,z'\to\infty$.
 
For bijections from $B_{k,l}^n$, $B_{k,\infty}^n$, $B_{\infty,\infty}^n$ we define a weakly increasing subsequence to be a sequence of pairs of associated elements, first element is from the first set, second element is from the second set, which weakly increase in each element. Under the RSK--correspondence the length of the longest weakly increasing subsequence of a bijection coincides with the length of the first row of the corresponding Young diagram in all three cases described above.   

\head \S3. Harmonic analysis on the infinite symmetric group
\endhead

For more detailed discussion of the material of this section see \cite{KOV}, \cite{VK}, \cite{P.I}.

We define the infinite symmetric group $S(\infty)$ as the inductive
limit of the finite symmetric groups $S_n$ with respect to natural
embeddings $S_n\to S_{n+1}$. Equivalently, $S(\infty)$ is the group of
{\it finite\/} permutations of the set $\{1,2,\dots\}$. 

By a {\it character\/} of $S(\infty)$ (in the sense of von Neumann)
we mean any central, positive definite function $\chi$ on
$S(\infty)$, normalized by the condition $\chi(e)=1$. We assign to
$\chi$ a function $M(\la)$ on the set $\Y=\sqcup \Y_n$ of Young
diagrams as follows: for any $n=1,2,\dots$,
$$
\chi\bigm|_{S_n}=\sum_{\la\in\Y_n} M(\la)\frac{\chi^\la}{\dim\la}\,,
$$
where $\chi^\la$ denotes the irreducible character of $S_n$ (in the
conventional sense), indexed by $\la\in\Y_n$, and
$\dim\la=\chi^\la(e)$ is its dimension. Let $M^{(n)}$ stand for
the restriction of the function $M$ to $\Y_n$; this is a probability
distribution on $\Y_n$. Conversely, let $M=\{M^{(n)}\}$ be a
function on $\Y$ such that each $M^{(n)}$ is a probability
distribution; then $M$ corresponds to a character $\chi$ if (and only
if) the distributions $M^{(n)}$ obey a natural coherence relation, which
comes from the classical Young branching rule for the irreducible
characters of the finite symmetric groups, see \cite{VK}, \cite{P.I}.
\footnote{Equivalently, the function
$\varphi(\la)=M(\la)/\dim\la$ must be a {\it harmonic function on the
Young graph\/} $\Y$ in the sense of Vershik and Kerov, see \cite{VK},
\cite{P.I}.}

\proclaim{Proposition 3.1} The z--measures $\Mnz$ introduced in \S1
satisfy the coherence relation mentioned above and, consequently,
define a character $\chi_{z,z'}$ of $S(\infty)$.
\endproclaim

Several direct proofs of the
proposition are known. E.g., a simple proof is given in \cite{P.I, \S7}. About generalizations, see \cite{K3}, \cite{BO3}. 

Note that the degenerations $M_{k,l}^{(n)}$, $M_{k,\infty}^{(n)}$, and
$M_{\infty,\infty}^{(n)}$ of the z--measures also correspond to certain
characters, which will be denoted as $\chi_{k,l}$, $\chi_{k,\infty}$,
and $\chi_{\infty,\infty}$, respectively. The character
$\chi_{\infty,\infty}$ is easily described: it takes value 1 at
$e\in S(\infty)$ and vanishes at all other elements of the group.

By the very definition of the characters of $S(\infty)$, they form a 
convex set. The extreme points of that set are called the {\it
indecomposable\/} characters, and the other points are called
{\it decomposable\/} characters. 

According to a remarkable theorem due to Thoma \cite{T1} (see also
\cite{VK}), the indecomposable characters of $S(\infty)$ are
parametrized by the points of the infinite dimensional simplex
$$
\Omega=\{\al_1\ge\al_2\ge\dots\ge 0,\ \be_1\ge\be_2\ge\dots\ge 
0\,|\,\sum_{i=1}^\infty(\al_i+\be_i)\le 1\},
$$
which is called the {\it Thoma simplex}. Given a point
$\omega=(\al,\be)\in\Omega$, we denote by $\chi^{(\omega)}$ the
corresponding indecomposable character.

The characters $\chi_{k,\infty}$ and $\chi_{\infty,\infty}$ are
indecomposable: the former corresponds to the point $\omega$ with 
$\al_1=\al_2=\dots=\al_k=1/k$ (all other coordinates are zero), and
the latter --- to the point $\omega=(0,0)$ (all coordinates are zero).  

The characters $\chi_{z,z'}$ (with $z,z'$ satisfying the conditions of \S1) and $\chi_{k,l}$ are decomposable.
 
Every character can be uniquely represented as a convex combination of the 
indecomposable ones,
$$
\chi=\int_\Omega \chi^{(\omega)}P(d\omega).
$$ 
Here $P$ is a probability measure on $\Omega$, which is called the
{\it spectral measure} of the character $\chi$. Moreover, any
probability measure on $\Omega$ is a spectral measure of a character,
so that the set of characters of $S(\infty)$ is isomorphic, as a
convex set, to the set of probability measures on the Thoma simplex.
Under this isomorphism, indecomposable characters correspond to delta
measures on $\Omega$. 

Given a concrete decomposable character $\chi$, a natural problem is
to describe explicitly its spectral measure $P$. This will be referred
to as the {\it problem of harmonic analysis}. 

This problem is readily solved for the degenerate characters
$\chi_{k,l}$: 

\proclaim{Proposition 3.2} Let $\chi=\chi_{k,l}$ with $k\le l$ and
set $a=l-k$. Then the spectral measure is concentrated on the
$(k-1)$-dimensional subsimplex 
$$
\al_1+\dots+\al_k=1, 
\qquad \al_{k+1}=\al_{k+2}=\dots=\be_1=\be_2\dots=0 
$$
of $\Omega$ and has density
$$
const\cdot \prod_{1\le i<j\le k}(\al_i-\al_j)^2 \cdot
\prod_{i=1}^k \al_i^a
$$
with respect to the Lebesgue measure.
\endproclaim

For the characters $\chi_{z,z'}$ with nonintegral parameters the
problem of harmonic analysis is highly nontrivial and will be
briefly discussed at the end of \S8. One of the first results in this direction is as
follows: 

\proclaim{Proposition 3.3} Let $P_{z,z'}$ denote the spectral measure of
$\chi_{z,z'}$. Except the obvious equality $P_{z,z'}=P_{z',z}$, the
measures $P_{z,z'}$ are pairwise disjoint. \footnote{Two measures are
called disjoint if there exist disjoint Borel sets supporting them.}
\endproclaim

Notice the following general result which relates the spectral
measure $P$ of a character $\chi$ to the finite probability 
distributions $M^{(n)}$. Let us embed $\Y_n$ into $\Omega$ by:
$$
\gathered
\la=(p_1,\dots,p_d\,|\,q_1,\dots,q_d)\in\Y_n\\ 
\mapsto\left\{\frac{p_1+1/2}n,\dots\frac{p_d+1/2}n,\,0,0,\dots;\, 
\frac{q_1+1/2}n,\dots,\frac{q_d+1/2}n, \,0,0,\dots\right\}\in\Omega.
\endgathered
$$

\proclaim{Proposition 3.4} As $n\to\infty$, the push--forwards of the
measures $M^{(n)}$ under these embeddings weakly converge to $P$. 
\endproclaim

This is a special case of a more general result proved in
\cite{KOO}. 

The characters of $S(\infty)$ can be related to representations in
two ways.

The first way is rather evident. Each character $\chi$ is a positive
definite function on $S(\infty)$, so that it determines a unitary
representation of $S(\infty)$, which will be denoted as $\Pi(\chi)$.
When $\chi$ is indecomposable, $\Pi(\chi)$ is a {\it factor\/}
representation of finite type in the sense of von Neumann, see
\cite{T2}.  

The second way is a bit more involved. Set $G=S(\infty)\times
S(\infty)$ and let $K$ denote the diagonal subgroup in $G$, which is
isomorphic to $S(\infty)$. We interpret $\chi$ as a function on the
first copy of $S(\infty)$, which is a subgroup of $G$, and then
extend it to the whole group $G$ by the formula
$$
\psi(g_1,g_2)=\chi(g_1g_2^{-1}), \qquad (g_1,g_2)\in G.
$$
Note that $\psi$ is the only extension of $\chi$ that is a
$K$-biinvariant function on $G$. The function $\psi$ is also positive
definite, so that one can assign to it a unitary representation in
the canonical way. This representation of the group $G$ will be
denoted by $T(\chi)$. By the very construction, it possesses a
distinguished $K$-invariant vector.

Note that $\Pi(\chi)$ coincides with the restriction of $T(\chi)$ to
the first copy of $S(\infty)$. If $\chi$ is indecomposable,
$\chi=\chi^{(\omega)}$, then $T(\chi)=T(\chi^{(\omega)})$ is
irreducible. The representations of the form $T(\chi^{(\omega)})$ are
exactly the irreducible unitary representations of the group $G$
possessing a $K$-invariant vector (such a vector is unique, within a
scalar factor). Thus, the Thoma simplex can be identified with the
{\it spherical dual\/} to $(G,K)$. \footnote{It is worth noting that the
irreducible representations of the form $T(\chi^{(\omega)})$ (except two
trivial cases) are {\it not\/} tensor products of irreducible
representations of the factors $S(\infty)$.} 

The representation $T(\chi_{\infty,\infty})$ is readily described: it
coincides with the natural representation of the group $G$
realized in the Hilbert space $\ell^2(G/K)$. Note that $G/K$ is
identified with the group $S(\infty)$ on which $G$ acts by left
and right shifts, so that $T(\chi_{\infty,\infty})$ may be called the
{\it regular\/} representation of $G$. As for
$\Pi(\chi_{\infty,\infty})$, it provides a classical realization of the
hyperfinite von Neumann factor of type $\text{II}_1$. 

The representations $T(\chi_{z,z'})$ are called the {\it generalized
regular\/} representations of $G$. The term is motivated by the fact
that each $T(\chi_{z,z'})$ can be realized as the inductive limit of
a chain of the form 
$$
\dots \to (\operatorname{Reg}_n, v_n)\to 
(\operatorname{Reg}_{n+1}, v_{n+1})\to\dots,
$$
where $\operatorname{Reg}_n$ stands for the (bi)regular
representation of the group $S_n\times S_n$ in the space of functions
on $S_n$ and $v_n$ is a certain vector in that space, depending on
the parameters $z,z'$ ($v_n$ is given by a certain central function
on $S_n$). When $z'=\bar z$, the generalized regular representations
admit a very nice realization --- in certain $L^2$ spaces of functions
defined on a compactification of the group $S(\infty)$. We refer to
\cite{KOV} for the exposition of this construction.

Finally, note that for any decomposable character $\chi$, the
decomposition of $T(\chi)$ into irreducible representations is
governed by the spectral measure $P$:
$$
T(\chi)=\int_\Omega T(\chi^{(\omega)})P(d\omega).
$$

\head \S4. Mixing \endhead

Theorem 1.2 computes the correlation functions of a point process obtained from the distributions $\Mnz$ mixed together by the negative binomial distribution with parameters $(t,\xi)$, see \S1. In this section we consider the degenerations of the mixing procedure in the cases when $z$ and $z'$ are integers, $z$ is an integer and $z'\to\infty$, $z$ and $z'$ both tend to infinity,
and $\xi\to 1$. 

If $z=k$ and $z'=l$ are positive integers then nothing interesting happens --- we have to mix the corresponding measures on $\Y_n$'s by the negative binomial distribution with parameters $(kl,\xi)$.

If $z=k$ is a positive integer and $z'\to\infty$, or $z\to\infty$ and $z'\to\infty$, $t=zz'$ goes to $\infty$. If we keep $\theta=t\xi$ fixed (hence, $\xi\to 0$) then the negative binomial distribution degenerates to the Poisson distribution with parameter $\theta$. The mixing procedure with Poisson distribution is called {\it poissonization}. 
 
The degeneration $\xi\to 1$ is a bit more delicate. 
Let us  embed $\Z'$ into the punctured line $\R^*=\R\setminus\{0\}$
and then rescale the lattice by multiplying the coordinates of
its points by $(1-\xi)$. Then the coordinates of the point configuration in $\R^*$ that corresponds to $\la\in \Y_n$  as defined in (1.2) after rescaling
differ from the coordinates of the image of $\la$ in $\Omega$ by the scaling
factor $(1-\xi)n$.  

The discrete distribution on the positive semiaxis with 
$$
\operatorname{Prob}\{(1-\xi)n\}=(1-\xi)^t\,\frac{(t)_n}{n!}\,\xi^n,\quad
n=0,1,2,\dots, 
$$ 
which depends on the parameter $\xi\in(0,1)$, converges, as
$\xi\to 1$, to the gamma--distribution with parameter $t$ 
$$
\gamma(ds)=\frac {s^{t-1}}{\Gamma(t)}\,e^{-s}ds.
$$

This brings us to the following construction.
Consider the space $\wt\Omega=\Omega\times \R_+$ with the probability
measure 
$$
\wt P_{z,z'}=P_{z,z'}\otimes \frac {s^{t-1}}{\Gamma(t)}\,e^{-s}ds.
$$
Let us embed $\Y=\Y_0\sqcup\Y_1\sqcup\Y_2\sqcup\dots$ into $\wt\Omega=\Omega\times \R_+$ by sending a Young diagram $\la\in\Y_n$ to the pair consisting of its image in $\Omega$ and the number $(1-\xi)n$. 
\proclaim{Proposition 4.1}
The push--forward of $\Mzxi$ under the embeddings described above converges, as $\xi\to 1$, to $\wt P_{z,z'}$.
\endproclaim
Exact claims with a detailed description of this convergence will appear in \cite{BO3}.

\head \S5. Ensembles \endhead
Before going further, let us introduce several terms which will be used below.

By the word {\it ensemble} throughout this paper we will mean a stochastic point process (i.e, a probability measure on the space of point configurations)  whose correlation functions
$\rho_n(x_1,\dots,x_n)$ are given by determinantal formulas of the form
$$
\rho_n(x_1,\dots,x_n)=\det[K(x_i,x_j)]_{i,j=1}^n
$$
where $K(x,y)$ is a certain kernel. We will call $K(x,y)$ the {\it correlation kernel}.

The process $\Pzxi$ is an example, the ensemble lives on $\Z'$  and the correlation kernel is the hypergeometric
kernel, see Theorem 1.2. We will call it the {\it discrete z--ensemble}. 

In all our examples the points of the ensembles will vary in  discrete or continuous subsets of the real line. Such a subset will be called the {\it phase space} of the corresponding ensemble. For example, $\Z'$ is the phase spase of  $\Pzxi$.

There is a class of {\it orthogonal polynomial ensembles} characterized by the condition of having a fixed finite number of points, say $k$, the joint probability distribution of which has the density
$$
const\cdot \prod_{1\le i,j\le k}(x_i-x_j)^2\prod_{i=1}^k w(x_i)
$$
with respect to either the Lebesgue measure if the phase space is continuous, or counting measure if the phase space is discrete. A standard argument due to Dyson \cite{Dy}, \cite{Me}, shows that the correlation kernel is the Christoffel--Darboux kernel of order $k$ for orthogonal polynomials on the phase space with respect to the weight function $w(x)$. If these polynomials have the form
$$
p_n(x)=a_nx^n+\{\text{lower degree terms}\}
$$
with $h_n=\Vert p_n\Vert^2$ then the kernel has the form
$$
K(x,y)=\frac {a_{k-1}}{a_kh_{k-1}}\,\frac{p_k(x)p_{k-1}(y)-p_{k-1}(x)p_k(y)}{x-y}\,\sqrt{w(x)w(y)}.
$$

Below we will consider the following orthogonal polynomial ensembles:

$\bullet$ {\it Laguerre ensemble}: phase space $\R_+$, weight function $w(x)=x^a e^{-x}$, $a>{-1}$;

$\bullet$ {\it Hermite ensemble}: phase space $\R$, weight function $w(x)=e^{-x^2}$;

$\bullet$ {\it Charlier ensemble}: phase space $\Z_+$, weight function $w(x)={\theta^x}/{x!}$, $\theta>0$;

$\bullet$ {\it Meixner ensemble}: phase space $\Z_+$, weight function $w(x)={(a+1)_x\xi^x}/{x!}$, $a>-1$, $\xi\in(0,1)$.

Corresponding normalizing constants for the orthogonal polynomials can be found in \cite{KS}, \cite{NSU}. The Christoffel--Darboux kernels for these ensembles will be called {\it Laguerre, Hermite, Charlier}, and {\it Meixner kernels}, respectively.

We will also deal with the {\it Airy ensemble}, see \cite{F1}, \cite{TW1}: the phase space is $\R$, the correlation kernel is
$$
\frac{A(x)A'(y)-A'(x)A(y)}{x-y}
$$
where $A(x)$ is the Airy function.

Two other ensembles that we will need are the ensemble arising from poissonized Plancherel distributions (see \cite{BOO} and \S6 below) with the phase space $\Z'$ and the kernel of the form (1.3), (1.4) where 
$$
P_\pm(x)=\theta^{\frac 14} J_x(2\sqrt{\theta}),\quad Q_\pm(x)=\theta^{\frac 14} J_{x+1}(2\sqrt{\theta}),
\tag 5.1
$$
$\theta>0$ is a parameter, $J_\nu(x)$ is the Bessel function; and the ensemble arising from the problem of harmonic analysis on $S(\infty)$ described in \S3 (see \cite{BO1} and \S6 below) with the phase space $\R^*$ and the kernel of the form (1.3) where
$$
\gathered 
F_1(x)=\cases-{\Cal Q_+(x)} ,&x>0\\  
{\Cal P_-(-x)},&x<0\endcases;\qquad
F_2(x)=\cases {\Cal P_+(x)} ,&x>0\\ 
{\Cal Q_-(-x)} ,&x<0\endcases;
\\
G_1(x)=\cases {\Cal P_+(x)},&x>0\\  
-{\Cal Q_-(-x)},&x<0\endcases;\qquad
G_2(x)=\cases {\Cal Q_+(x)} ,&x>0\\  
{\Cal P_-(-x)},&x<0\endcases,
\endgathered
$$
$$
\gathered
\Cal P_\pm(x)=\frac{(zz')^{1/4}}{(\Ga(1\pm z)\Ga(1\pm z')\,x)^{1/2}}\,
W_{\frac{\pm(z+z')+1}2,\frac{z-z'}2}(x), \\
\Cal Q_\pm(x)=\frac{(zz')^{3/4}}{(\Ga(1\pm z)\Ga(1\pm z')\,x)^{1/2}}\,
W_{\frac{\pm(z+z')-1}2,\frac{z-z'}2}(x),
\endgathered  
$$  
$z,z'$ satisfy the assumptions stated in \S1, $W_{\kappa,\mu}(x)$ is the Whittaker function. 
We will call these ensembles the {\it Plancherel ensemble} and the {\it continuous z--ensemble} respectively. The kernel for the first one will be called the {\it Plancherel kernel}, for the second one --- the {\it Whittaker kernel}.
  
When an ensemble lives on $\R^*$ or $\Z'$, one may single out
its {\it positive part\/} --- the restriction to $\R_+\subset\R^*$ or $\Z'_+\subset\Z'$, respectively.
The correlation kernel of the positive part is the corresponding restriction of the correlation kernel of the initial ensemble. We will use the term ``positive part of the kernel'' for such restrictions. 

The positive part of the Plancherel kernel has been independently found in \cite{J2} where it was called the {\it discrete Bessel kernel}, see \S9.

\head \S6. Correlations after mixing \endhead

In view of \S2, it is natural to denote the measures on Young diagrams with $n$ boxes coming from $B_{k,l}^n$, $B_{k,\infty}^n$, $B_{\infty,\infty}^n$ as $M_{k,l}^{(n)}$, $M_{k,\infty}^{(n)}$, $M_{\infty,\infty}^{(n)}$, respectively, and the corresponding mixtures (i.e., measures on the set of all Young diagrams) as $M_{k,l,\xi}$, $M_{k,\infty,\theta}$, $M_{\infty,\infty,\theta}$.
We want to see how the hypergeometric kernel will behave in these degenerate cases.

Let us start with the case when $z,z'$ are positive integers,
say, $z=k$, $z'=l$, $k\le l$. Denote $a=l-k$.
  
\proclaim{Proposition 6.1 (\cite{BO2}, \cite{J1})} Let $\la=(p_1,\dots,p_d\,|\,q_1,\dots,q_d)\in \Y$ be distributed according to $M_{k,l,\xi}$. Then the distribution of points $\{k+p_1,\dots,k+p_d\}$ coincides with the restriction of the Meixner ensemble with parameters $(a,\xi)$ to the set $\{k,k+1,\dots\}$.
\endproclaim 
 
This claim corresponds to the fact that the hypergeometric functions participating in the hypergeometric kernel become Meixner polynomials if $z$ or $z'$ is integral, see \cite{BO2}. Furthermore, the positive part of the hypergeometric kernel becomes the Christoffel--Darboux kernel for Meixner polynomials (shifted by $k$). 

Now we pass to $M_{k,\infty,\theta}$.

\proclaim{Proposition 6.2 (\cite{J2})} Let $\la=(p_1,\dots,p_d\,|\,q_1,\dots,q_d)\in \Y$ be distributed according to $M_{k,\infty,\theta}$. Then the distribution of points $\{k+p_1,\dots,k+p_d\}$ coincides with the restriction of the Charlier ensemble with parameter $\theta$ to the set $\{k,k+1,\dots\}$.
\endproclaim

The easiest way to see this is to observe the degeneration of Meixner polynomials with parameters $(a,\xi)$ to Charlier polynomials with parameter $\theta$ when $a\to\infty$, $\theta=k(k+a)\xi$ is fixed.

Next, consider the situation when $z$ and $z'$ both go to $\infty$.
\proclaim{Proposition 6.3 (\cite{BOO})} Let $\la=(p_1,\dots,p_d\,|\,q_1,\dots,q_d)\in \Y$ be distributed according to $M_{\infty,\infty,\theta}$. Then the random point configuration $\{p_1+\frac 12,\dots,p_d+\frac 12,-q_1-\frac 12,\dots,-q_d-\frac 12\}$ forms the Plancherel ensemble with parameter $\theta$.
\endproclaim

This claim corresponds to the degeneration of the hypergeometric function to the Bessel $J$--function when first two parameters go to infinity and the argument goes to zero so that the product of these three numbers is fixed (and equals $\theta$). 

As for the representation theoretic picture, we have the following claim.

\proclaim{Proposition 6.4 (\cite{BO1})} Let $((\al,\be),s)\in \wt\Omega=\Omega\times \R_+$ be distributed according to $\wt P_{z,z'}$. Then the random point configuration $(s\al_1,s\al_2,\dots,-s\be_1,-s\be_2,\dots)$ forms the continuous z--ensemble.
\endproclaim

\example{Remark 6.5} When one of the parameters $z,z'$ becomes integral, say, $z=k\in\{1,2,\dots\}$, and $z'=z+a$, $a>-1$, the Whittaker kernel degenerates to the Laguerre kernel of order $k$ with parameter $a$. Then Proposition 6.4 implies that the measure $\wt P_{z,z'}$ gets concentrated on the finite--dimensional subset of $\wt\Omega=\Omega\times \R_+$ where $\al_{k+1}=\al_{k+2}=\dots=\be_1=\be_2=\dots=0$, and on this subset in the new coordinates $x_i=s\al_i$ ($s$ is the coordinate on $\R_+$) it equals, see \cite{P.III, Remark 2.4},
$$
const\cdot \prod_{1\le i<j\le k}(x_i-x_j)^2\prod_{i=1}^kx_i^{a}e^{-x_i}dx_i.
$$
This agrees with Proposition 3.2.
\endexample

\head \S7. Asymptotics when mixing parameters tend to a limit \endhead

We start with $M_{k,l,\xi}$. Assume that $a=l-k\ge 0$. Then Proposition 3.4, the degeneration of the hypergeometric kernel to the Whittaker kernel (Proposition 6.4) and the coincidence of the Whittaker kernel with the Laguerre kernel when at least one parameter is integral (Remark 6.5) justify the following claim.

\proclaim{Proposition 7.1} Let $\la\in \Y$ be distributed according to $M_{k,l,\xi}$. Then the random point configuration $\{(1-\xi)\la_1,\dots,(1-\xi)\la_k\}$ converges, as $\xi\to 1$, to the Laguerre ensemble.
\endproclaim

Now let us pass to $M_{k,\infty,\theta}$. The fact that the character of $S(\infty)$ corresponding to $M_{k,\infty}^{(n)}$ is indecomposable and corresponds to the point $\al_1=\dots=\al_k=1/k$ in $\Omega$ (see \S3) leads to the following statement.

Consider the embedding of the set of Young diagrams with length $\le k$ into $\R_+^k$ defined by normalizing the lengths of rows of a Young diagram by $\theta$. 

\proclaim{Proposition 7.2} Under the embeddings described above
$M_{k,\infty,\theta}$ weakly converges to the delta measure at the point $(1/k,\dots,1/k)$ as $\theta\to\infty$.
\endproclaim

One can also ask about fluctuations of $M_{k,\infty,\theta}$ around the limit delta measure.
Johansson \cite{J2} proved the following statement.

\proclaim{Proposition 7.3 (\cite{J2})} Let $\la\in\Y$ be distributed according to $M_{k,\infty,\theta}$. Then the random point configuration 
$$
\left\{\frac{\la_1-\theta/k}{\sqrt{2\theta/k}},\dots,\frac{\la_k-\theta/k}{\sqrt{2\theta/k}}\right\}
$$
converges, as $\theta\to\infty$, to the Hermite ensemble.
\endproclaim

The convergence of distribution of the first point of the random configuration from Proposition 7.3 was proved by Tracy and Widom, \cite{TW4}.

Propositions 7.2 and 7.3 correspond to a certain degeneration of Charlier polynomials to Hermite polynomials which follows from a more general degeneration of Laguerre polynomials with large argument and parameter to Hermite polynomials, see \cite{Te}. 

The most interesting case is $M_{\infty,\infty,\theta}$. The reason is simple -- the number of points (rows of Young diagrams) is not bounded in this case. 
One can look at at least two different regimes when $\theta\to\infty$: ``in the bulk of spectrum'' or ``at the edge of spectrum''. 

\proclaim{Proposition 7.4 (\cite{BOO}, \cite{J2})} Let $\la\in \Y$ be distributed according to $M_{\infty,\infty,\theta}$. Then the random point configuration
$$
\left\{\frac{\la_1-2\sqrt{\theta}}{\theta^{\frac 16}},\frac{\la_2-2\sqrt{\theta}}{\theta^{\frac 16}},\dots\right\}
$$
converges, as $\theta\to+\infty$, to the Airy ensemble.
\endproclaim

The convergence of distributions of the first and the second points of the random configuration from Proposition 7.4 was proved earlier in \cite{BDJ1}, \cite{BDJ2}.

Proposition 7.4 is the result of degeneration of the Bessel functions (5.1) to the Airy function and its derivative.  

For the results on the asymptotics ``in the bulk of spectrum''  we refer the reader to \cite{BOO}. These results correspond to the degeneration of the Plancherel kernel to the {\it discrete sine kernel} 
$$
\frac{\sin(a(x-y))}{\pi(x-y)},\qquad x,y\in \Z,\quad 0<a<\pi.
$$

It was also mentioned in \cite{J2} that under a certain limit procedure the Plancherel kernel degenerates to the conventional sine kernel $\sin(\pi(x-y))/(\pi(x-y))$ on $\R$. 

One can also consider ``double limits'' of $M_{k,l,\xi}$ and $M_{k,\infty,\theta}$ (or, equivalently, Meixner and Charlier ensembles) when at least two parameters tend to a critical value. Then the scaling procedure must involve at least two large parameters. For $M_{k,\infty,\theta}$ the asymptotics looks as follows.

\proclaim{Proposition 7.5 (\cite{J2})} Let $\la\in\Y$ be distributed according to $M_{k,\infty,\theta}$. Then the random point configuration
$$
\left\{\frac{\la_1-\theta/k-2\sqrt{\theta}}{(1+\sqrt{\theta}/k)^{\frac 23}\theta^{\frac 16}},\frac{\la_2-\theta/k-2\sqrt{\theta}}{(1+\sqrt{\theta}/k)^{\frac 23}\theta^{\frac 16}},\dots\right\}
$$
converges, as $k\to\infty$ and $\theta\to\infty$, to the Airy ensemble.
\endproclaim
 
The result corresponds to a degeneration of Charlier polynomials to the Airy function \cite{J2}.   
 
For $M_{k,l,\xi}$ a similar result was proved for $k,l\to+\infty$ in \cite{J1}. 
 
\head \S8. Asymptotics of non--mixed measures for large $n$
\endhead

As we have seen above, after mixing the study of our measures is not very difficult --- we just need to look at the corresponding degenerations of the hypergeometric kernel.
The picture before mixing is more subtle. 

For $M_{k,l}^{(n)}$ and $M_{k,\infty}^{(n)}$ the asymptotics before and after mixing are different. In comparison to the mixed cases, there appear restrictions on the supports of the limit measures. These restrictions come from the trivial condition that the sum of lengths of rows of a Young diagram with $n$ boxes is equal to $n$.

Consider the embedding of the set of Young diagrams with $n$ boxes and length $\le k$ into $\R_+^k$ defined by normalizing the lengths of rows of a Young diagram by $n$.

Proposition 3.2 and Proposition 3.4 lead to the following claim.

\proclaim{Proposition 8.1} As $n\to\infty$, the images of the measures $M_{k,l}^{(n)}$ under the embeddings
defined above converge to a measure concentrated on the set 
$$
\{(x_1,\dots,x_k)\in\R_+^k\,|\,x_1\ge x_2\ge\dots\ge x_k,\ \sum_{i=1}^kx_i=1\}.
$$
The density of the limit
measure with respect to the Lebesgue measure equals
$$
const\cdot \prod_{1\le i<j\le k}(x_i-x_j)^2\prod_{i=1}^kx_i^{a}
$$
(recall that $a=l-k\ge 0$).
\endproclaim

Similarly to Proposition 7.2, we have
\proclaim{Proposition 8.2} As $n\to\infty$, the images of the measures $M_{k,\infty}^{(n)}$ under the embeddings
defined above converge to the delta measure at the point $(1/k,\dots,1/k)$.
\endproclaim

Again, the fluctuations around the limit delta measure were determined by Johansson \cite{J2}. 

Define an embedding of the set of Young diagrams with $n$ boxes and length $\le k$ into $\R^k$ setting the $i$th coordinate of the image of $\la\in \Y_n$ equal to $\frac{\la_i-n/k}{\sqrt{2n/k}}$, cf. Proposition 7.3.

\proclaim{Proposition 8.3 (\cite{J2})}As $n\to\infty$, the images of the measures $M_{k,\infty}^{(n)}$ under the embeddings
defined above converge to a measure concentrated on the set 
$$
\{(x_1,\dots,x_k)\in\R^k\,|\,x_1\ge x_2\ge\dots\ge x_k,\ \sum_{i=1}^k x_i=0\}.
$$
The density of the limit
measure with respect to the Lebesgue measure equals
$$
const\cdot \prod_{1\le i<j\le k}(x_i-x_j)^2\cdot e^{-x_1^2-\dots-x_k^2}.
$$
\endproclaim

For the values of $M_{k,\infty}^{(n)}$ on functions depending only on $\la_1$ the claim was proved by Tracy and Widom \cite{TW4}.

In a sense, $M_{\infty,\infty}^{(n)}$ is the most pleasant measure. In this case the asymptotics of $M_{\infty,\infty}^{(n)}$ in the bulk of spectrum and at the edge of spectrum as $n\to\infty$ is exactly the same as the asymptotics of $M_{\infty,\infty,\theta}$ as $\theta\to\infty$.  We can say that the asymptotics admits {\it depoissonization}, see \cite{BOO} and also \cite{J2}. Let us explicitly state the analog of Proposition 7.4.

\proclaim{Proposition 8.4 (\cite{BOO}, \cite{J2})}
Let $\la\in \Y_n$ be distributed according to $M_{\infty,\infty}^{(n)}$. Then the random point configuration
$$
\left\{\frac{\la_1-2\sqrt{n}}{n^{\frac 16}},\frac{\la_2-2\sqrt{n}}{n^{\frac 16}},\dots\right\}
$$
converges, as $n\to+\infty$, to the Airy ensemble.
\endproclaim

Again, the convergence of distributions of first two points was proved in \cite{BDJ1}, \cite{BDJ2}.

Depoissonization of the result in the bulk of spectrum requires different ideas from those used in the proof of Proposition 8.4.
For the discussion of this case we refer to \cite{BOO}.

Proposition 7.5 also admits depoissonization.

\proclaim{Proposition 8.5 (\cite{J2})} Let $\la\in\Y_n$ be distributed according to $M_{k,\infty}^{(n)}$. Then the random point configuration
$$
\left\{\frac{\la_1-n/k-2\sqrt{n}}{(1+\sqrt{n}/k)^{\frac 23}n^{\frac 16}},\frac{\la_2-n/k-2\sqrt{n}}{(1+\sqrt{n}/k)^{\frac 23}n^{\frac 16}},\dots\right\}
$$
converges, as $k\to\infty$, $n\to\infty$ so that $(\ln n)^{\frac 16}/k\to 0$, to the Airy ensemble.
\endproclaim

The structure of spectral z--measures $P_{z,z'}$ defined in \S3 for general $z$ and $z'$ is fairly complicated. Note that $P_{z,z'}$ is the limit of the $n$th level z--measures $M_{z,z'}^{(n)}$, see Proposition 3.4.

Every probability measure on $\Omega$ (definition in \S3) can be viewed as a point process on $\R^*$, if we associate to every point $(\al,\be)\in\Omega$ the point configuration
$(\al_1,\al_2,\dots,-\be_1,-\be_2,\dots)$, cf. Proposition 6.4. The correlation functions of the process corresponding to $P_{z,z'}$ were all explicitly computed in \cite{P.II}. They do not have determinantal form and can be expressed through multivariate hypergeometric functions.

The situation after mixing is substantially simpler: the process associated to $\wt P_{z,z'}$ is the Whittaker ensemble (Proposition 6.4).

We refer to \cite{BO1}, \cite{P.I--V} for detailed discussion of measures $P_{z,z'}$, $\wt P_{z,z'}$ and associated point processes. 
 
\head \S9. Limit transitions \endhead 

The fact that numerous kernels and ensembles described above originated from the same hypergeometric kernel suggests a number of different limit transitions between them. 

On the top of the hierarchy we have the hypergeometric kernel which degenerates to all ensembles described above. This corresponds to the fact that the hypergeometric function is on the top of the hierarchy of classical special functions in one variable.
The kernel depends on three parameters $z,z',\xi$, and lives on the lattice $\Z'$.

The Meixner kernel is the specialization of the positive part of the hypergeometric kernel when one of the parameters $z,z'$ is integral. To be concrete, we will assume below that $z\in \{1,2,\dots\}$.

The Charlier kernel and the Whittaker kernel are one step below --- they both depend on two parameters, $(z,\theta)$ and $(z,z')$, respectively. The Charlier kernel is obtained from the Meixner kernel by taking the limit $z'-z\to+\infty$ with $\theta=zz'\xi$ fixed, the Whittaker kernel is obtained from the hypergeometric kernel via a scaling limit when $\xi\to 1$. The Charlier kernel lives on $\Z_+$, the Whittaker kernel lives on $\R^*$.

The Laguerre kernel is a particular case of the positive part of the Whittaker kernel when one of parameters $(z,z')$ is integral. It can be also obtained from the Meixner kernel by taking the limit $\xi\to 1$ (Proposition 7.1). The Laguerre kernel depends on two parameters $(z,a=z'-z)$ and lives on $\R_+$.

The Plancherel kernel is on the next level --- it lives on $\Z'$, depends on one parameter $\theta$ and can be obtained from the hypergeometric kernel via the limit $z,z'\to\infty$, $\xi\to 0$, $\theta=zz'\xi$ fixed. Its positive part can be obtained either from the Meixner kernel by letting $z,z'\to\infty$ with $\theta=zz'\xi$ fixed, or from the Charlier kernel by taking the limit $z\to\infty$. These two transitions are thoroughly discussed in \cite{J2}. 
It is worth noting that the whole Plancherel ensemble, as opposed to its positive part, cannot be obtained by taking limits of Meixner or Charlier ensembles.    

The Hermite kernel also depending on one integral parameter $z$ can be obtained from the Charlier kernel via the limit $\theta\to\infty$ (Proposition 7.3).  

The Airy kernel is at the bottom --- it has no parameters. It can be obtained in a number of different ways. For example, one can obtain the Airy kernel in the limit $\theta\to+\infty$ of the Plancherel kernel at the edge of spectrum (Proposition 7.4), or as the limit at the edge of spectrum of the Hermite kernel and the Laguerre kernel with parameter $a$ fixed when the order $z$ of these polynomial ensembles goes to infinity, \cite{F1}, \cite{TW1}. It can also be obtained as a double limit of Charlier or Meixner kernels, see the end of \S7, \cite{J1}, \cite{J2}.  

Of course, this is not the end of the story. The discrete sine kernel and the conventional sine kernel can be obtained from Plancherel kernel as $\theta\to\infty$, see \S7. The so--called Bessel kernel can be extracted from the Laguerre kernel ``at the hard edge of spectrum'' \cite{F1}, \cite{NW2}, \cite{TW2}. The sine kernel can be obtained from the Laguerre and Hermite kernels in the bulk of spectrum, see, e.g., \cite{NW1}. A number of new kernels can be obtained from the Whittaker kernel, see \cite{P.V}. Presumably, all these kernels can also be obtained as double or triple limits of the hypergeometric kernel.

Thus, a variety of kernels known so far can be obtained from the hypergeometric kernel, often in several different ways. As we tried to demonstrate above, sometimes such degenerations also carry the information about the asymptotic behavior of certain combinatorial objects.

\bigskip
\bigskip

{\smc A.~Borodin}: Department of Mathematics, The University of
Pennsylvania, Philadelphia, PA 19104-6395, U.S.A.  

E-mail address:
{\tt borodine\@math.upenn.edu}

{\smc G.~Olshanski}: Dobrushin Mathematics Laboratory, Institute for
Problems of Information Transmission, Bolshoy Karetny 19, 101447
Moscow GSP-4, RUSSIA.  

E-mail address: {\tt olsh\@iitp.ru,
olsh\@glasnet.ru}

\enddocument